\documentclass[english,pre,preprint,noshowpacs]{revtex4}
\usepackage[T1]{fontenc}
\usepackage[latin9]{inputenc}
\usepackage{units}
\usepackage{amsmath}
\usepackage{amssymb}
\usepackage{graphicx}
\usepackage{esint}

\makeatletter
\@ifundefined{textcolor}{}
{%
 \definecolor{BLACK}{gray}{0}
 \definecolor{WHITE}{gray}{1}
 \definecolor{RED}{rgb}{1,0,0}
 \definecolor{GREEN}{rgb}{0,1,0}
 \definecolor{BLUE}{rgb}{0,0,1}
 \definecolor{CYAN}{cmyk}{1,0,0,0}
 \definecolor{MAGENTA}{cmyk}{0,1,0,0}
 \definecolor{YELLOW}{cmyk}{0,0,1,0}
}

\usepackage{units}

\@ifundefined{textcolor}{}{%
 \definecolor{BLACK}{gray}{0}
 \definecolor{WHITE}{gray}{1}
 \definecolor{RED}{rgb}{1,0,0}
 \definecolor{GREEN}{rgb}{0,1,0}
 \definecolor{BLUE}{rgb}{0,0,1}
 \definecolor{CYAN}{cmyk}{1,0,0,0}
 \definecolor{MAGENTA}{cmyk}{0,1,0,0}
 \definecolor{YELLOW}{cmyk}{0,0,1,0}
}
\usepackage{enumitem}% customizable list environments
\usepackage[english]{babel}

\makeatother

\usepackage{babel}
\begin{document}

\title{Impact Mechanics of Elastic Structures with Point Contact}

\author{Róbert Szalai}

\date{January 18, 2014}

\affiliation{Department of Engineering Mathematics, University of Bristol, Queen's
Bldg., University Walk, Bristol, BS8 1TR, UK, Email: r.szalai@bristol.ac.uk}

\begin{abstract}
This paper introduces a modeling framework that is suitable to resolve
singularities of impact phenomena encountered in applications. The
method involves an exact transformation that turns the continuum,
often partial differential equation description of the contact problem
into a delay differential equation. The new form of the physical model
highlights the source of singularities and suggests a simple criterion
for regularity. To contrast singular and regular behavior the impacting
Euler-Bernoulli and Timoshenko beam models are compared. 
\end{abstract}
\maketitle

\section{Introduction}

Impact mechanics is a great concern for engineers, thus many models
were developed to understand this phenomena \cite{Stronge}. Even
simple models of impact \cite{KhuliefShabana1987} lead to complicated
predictions, like infinite chatter \cite{NordmarkPiiroinen2009},
period adding bifurcations, chaos \cite{ChinOttNusse1994} and non-deterministic
motion \cite{PainleveAlan}. The more elaborate models, however can
suffer from convergence problems as either the time step of the solution
decreases \cite{WangKim1996} or the resolved degrees of freedom is
increased \cite{Melcher}. This signals the need for a better modeling
framework for impact phenomena.

Most state-of-the-art impact models are finite dimensional and predict
infinite contact forces. This happens because at impact the contact
point abruptly changes its velocity and has an infinite acceleration.
If the mass of the contact point is not zero, this infinite acceleration
requires an infinite contact force. One example is the standard modal
description \cite{EwinsBook} of linear structures. Each vibration
mode is associated with a non-zero modal mass. Truncation of the modal
expansion describes the motion as if a finite set of rigid bodies
were coupled with springs and dampers. The arising infinite forces
are difficult to handle and are the source of singularities. In contrast,
contact points of elastic bodies have infinitesimally small mass.
This means that the contact force can stay finite despite infinite
accelerations. In this paper we show that finite contact forces are
possible under general conditions.

In some models, impact is treated as a discrete-time event. Since
impacting bodies must not overlap, the velocity state of the bodies
must be altered at contact, so that their subsequent motion avoids
overlap momentarily. One such change of velocities can described by
the coefficient of restitution (CoR) model that stipulates that the
incident velocity and a rebound velocity of the contact points are
opposite and linearly related. This still leaves a great deal freedom
in choosing the rest of the rebound velocities of the structure, therefore
the CoR model must be coupled with additional rules, which can be
based on momentum balances \cite{PalasShabana1992,Vyasarayani2010}
or collocation \cite{WaggBishop2002}. The main weakness of the CoR
method is that it leads to high-frequency chatter that is proportional
to the highest natural frequency of the system. This phenomenon restricts
the highest natural frequencies that can be included in the model,
because numerical simulation would becomes prohibitively slow. There
are a few extensions to the CoR method that avoid high-frequency chatter
\cite{NordmarkPiiroinen2009,Segalman2006}, however they do not fix
the cause of the problem which is the presence of infinite or undefined
contact forces. To illustrate this point we compare our method to
a CoR model \cite{Vyasarayani2010} and show how chatter is eliminated
when the contact force is well-defined.

Impact can also be modeled by calculating the contact force using
impulse response functions of the elastic structures at the contact
point. This approach is expected to be more accurate, however convergence
can be as troublesome as for CoR models. For example, Wang and Kim
\cite{WangKim1996} found that in the limit of zero time-step of their
algorithm, the contact force diverges when an Euler-Bernoulli beam
impacts an elastic St Venant rod. This divergence is related to the
approximation of the impulse-response function with finitely many
vibration modes. Similar approaches were used by other authors \cite{Evans1991,Fathi1994,YinQinZou2007}
using time stepping algorithms and finitely many vibration modes;
in some cases the shape of the contact force was assumed \cite{Langley2012},
which avoids divergence.

In this paper we use a similar technique to impulse-response functions,
and transform the governing equation of the impacting mechanical system
to a delay differential equation. The memory term of our equation
can be thought of as the convolution integral with the impulse response
function. Depending on the properties of the memory term the model
is either singular or regular. In case of a regular model the finite
contact force is uniquely calculated from a delay-differential equation.

\section{Mechanical model}

We analyze impact mechanics through a converging series expansion
of the continuum problem. We use infinitely many vibration modes $x_{k}(t)$
to recover the entire motion of the structure. Through this expansion
the displacement of any point $\chi$ of the elastic body can be written
as an infinite sum 
\[
u(t,\chi)=\sum_{k=1}^{\infty}\psi_{k}(\chi)x_{k}(t),
\]
where $\psi_{k}(\chi)$ are the mode shapes of the structure \cite{EwinsBook}.
If the motion of the structure is decoupled into non-resonant modes
of vibration, the equation of motion can be written as

\begin{equation}
\ddot{\boldsymbol{x}}(t)+2\boldsymbol{D}\boldsymbol{\Omega}\dot{\boldsymbol{x}}(t)+\boldsymbol{\Omega}^{2}\boldsymbol{x}(t)=\boldsymbol{f}_{e}(t)+\boldsymbol{n}f_{c}(t),\label{eq:modeDecomp}
\end{equation}
where $\boldsymbol{x}=(x_{1},x_{2},\ldots)^{T}$, the mass matrix
is assumed to be the identity, $\boldsymbol{\Omega}=\mathrm{diag}(\omega_{1},\omega_{2},\ldots)$
and $\boldsymbol{D}=\mathrm{diag}(D_{1,}D_{1,}\ldots)$, $\boldsymbol{f}_{e}(t)$
represents the external force, and $f_{c}(t)$ is the contact force.
We assume that the natural frequencies scale according to $\omega_{k}=\omega_{1}k^{\alpha}$,
for $k\gg1$. Vector $\boldsymbol{n}$ in equation (\ref{eq:modeDecomp})
represents the contribution of the modes to the motion of the contact
point $y(t)=\boldsymbol{n}\cdot\boldsymbol{x}(t)$ with 
\begin{equation}
\boldsymbol{n}=(\psi_{1}(\chi^{\star}),\psi_{2}(\chi^{\star}),\ldots)^{T},\label{eq:1stProj}
\end{equation}
where $\chi^{\star}$ represents the contact point. The method described
in this paper is not restricted to modal equations \eqref{eq:modeDecomp},
a more general description can be found in \cite{SzalaiMZfriction}.

\subsection{Approximating the contact force}

To better understand the impact process we first approximate the contact
force assuming that the impact is infinitesimally short. We assume
a single structure that interacts with a rigid stop. Contact occurs
at $t_{0}$ if $y(t_{0})=0$. As a first step we calculate a constant
contact force that keeps the stop penetrating the structure after
time $\delta t$, that is, $y(t_{0}+\delta t)=0$, which forms a boundary
value problem. After solving equation (\ref{eq:modeDecomp}) the contact
force becomes 
\[
f_{c}=-C\delta t^{\frac{1}{\alpha}-1}\left(\boldsymbol{n}\cdot\dot{\boldsymbol{x}}(t_{0}^{-})\right),
\]
where $\dot{\boldsymbol{x}}(t_{0}^{-})$ is the vector of modal velocities
just before the impact and $0<C<\infty$ is a constant. When $\delta t\to0$
the velocity of the contact point reverses and that corresponds to
a unit CoR. The details of the calculation can be found in the Appendix.

When evaluating the contact force there are three cases as $\delta t\to0$
at the onset of contact. If $\alpha<1$ the contact force becomes
zero, if $\alpha=1$, the contact force is a finite constant and for
$\alpha>1$ the contact force tends to infinity. This simple result
implies that for a finite contact force at least a subsequence of
the natural frequencies $\omega_{k}$ must scale at most linearly
as $k\to\infty$.

For a system composed of an elastic body which strikes a rigid stop,
the impact should change the momentum of the elastic body. Our calculation
shows that the change of momentum of the elastic body is zero, i.e.,
$\lim_{\delta t\to0}\left(f_{c}\delta t\right)=0$ as the contact
time tends to zero. This implies that the impact must occur during
a non-zero and finitely long time-interval.

\section{Model transformation}

To accurately calculate the contact force as a function of time in
the continuum problem we transform the infinite dimensional system
(\ref{eq:modeDecomp}) into a delay equation. Delay terms can naturally
arise from traveling wave solutions of partial differential equations
\cite{ShimmyDelay,StepanForgeDETC}. However dispersion might prevent
one to write down such solutions. Instead we use the Mori-Zwanzig
formalism as is described for mechanical systems in \cite{SzalaiMZfriction}
and obtain a time-delay model.

Our aim is to find a self-contained equation that exactly describes
the evolution of $\boldsymbol{y}(t)=\left(\boldsymbol{n}\cdot\boldsymbol{x}(t),\boldsymbol{n}\cdot\dot{\boldsymbol{x}}(t)\right)^{T}$.
We call $\boldsymbol{y}$ the vector of resolved variables. The first
step in the process is to transform (\ref{eq:modeDecomp}) into a
first-order form 
\begin{equation}
\dot{\boldsymbol{z}}(t)=\boldsymbol{R}\boldsymbol{z}(t)+\left(\begin{array}{c}
0\\
\boldsymbol{n}f_{c}(t)
\end{array}\right)+\boldsymbol{f}(t),\;\boldsymbol{R}=\left(\begin{array}{cc}
\boldsymbol{0} & \boldsymbol{I}\\
-\boldsymbol{\Omega}^{2} & -2\boldsymbol{D}\boldsymbol{\Omega}
\end{array}\right),\label{eq:1stOrder}
\end{equation}
where $\boldsymbol{f}(t)=(\boldsymbol{0},\boldsymbol{f}_{e}(t))^{T}$.
$ $ Note that $\boldsymbol{R}$ can also represent any convergent
expansion of the continuum problem, including numerical schemes such
as finite difference methods and $\boldsymbol{y}$ can include any
finite number of variables \cite{SzalaiMZfriction}. To arrive at
a model that describes the evolution of the resolved variables $\boldsymbol{y}$,
we construct a projection with a finite dimensional range with the
help of the matrices 
\begin{gather*}
\boldsymbol{V}=\left(\begin{array}{cc}
\boldsymbol{n}^{T} & 0\\
0 & \boldsymbol{n}^{T}
\end{array}\right)\;\mbox{and}\;\boldsymbol{W}=\left(\begin{array}{cc}
\boldsymbol{m} & 0\\
0 & \boldsymbol{m}
\end{array}\right),
\end{gather*}
where the vector $\boldsymbol{m}$ is chosen such that $\boldsymbol{m}\cdot\boldsymbol{n}=1$
and its components obey $\left[\boldsymbol{m}\right]_{j}=0$ for $j>M<\infty$.
To simplify our analysis we also assume that the columns of $\boldsymbol{W}$
span an invariant subspace of $\boldsymbol{R}$. The resolved variables
now can be expressed as $\boldsymbol{y}=\boldsymbol{V}\boldsymbol{z}$,
and the projection and the reciprocal projection matrices become $\boldsymbol{S}=\boldsymbol{W}\boldsymbol{V}$
and $\boldsymbol{Q}=\boldsymbol{I}-\boldsymbol{S}$, respectively.
Further, we assume that the initial condition of (\ref{eq:1stOrder})
is specified at $t=0$ and that there is no contact force initially.
According to \cite{SzalaiMZfriction}, with this notation the governing
equation for the resolved variables becomes 
\begin{equation}
\frac{\mathrm{d}}{\mathrm{d}t}\boldsymbol{y}(t)=\boldsymbol{A}\boldsymbol{y}(t)+\boldsymbol{L}^{\infty}f_{c}(t)+\int_{0}^{t}\mathrm{d}_{\tau}\boldsymbol{L}(\tau)\frac{\mathrm{d}}{\mathrm{d}t}\left[f_{c}(t-\tau)\right]+\boldsymbol{g}(t),\label{eq:redLinForce}
\end{equation}
where $\boldsymbol{A}=\boldsymbol{V}\boldsymbol{R}\boldsymbol{W}$,
the memory kernel is a function of bounded variation, 
\begin{align*}
\boldsymbol{L}(\tau) & =\int_{0}^{\tau}\left(\boldsymbol{V}\mathrm{e}^{\boldsymbol{R}\boldsymbol{Q}\theta}(0,\boldsymbol{n})^{T}-\boldsymbol{L}^{\infty}\right)\mathrm{d}\theta,\\
\boldsymbol{L}^{\infty} & =\boldsymbol{A}\boldsymbol{V}\boldsymbol{R}^{-1}(0,\boldsymbol{n})^{T},
\end{align*}
and the forcing term is 
\[
\boldsymbol{g}(t)=\boldsymbol{V}\boldsymbol{R}\boldsymbol{Q}\mathrm{e}^{\boldsymbol{R}t}(\boldsymbol{x}(0),\dot{\boldsymbol{x}}(0))^{T}+\int_{0}^{t}\left(\boldsymbol{V}\boldsymbol{R}-\boldsymbol{A}\boldsymbol{V}\right)\mathrm{e}^{\boldsymbol{R}\tau}\boldsymbol{f}(t-\tau)\mathrm{d}\tau.
\]
The integral in (\ref{eq:redLinForce}) is meant in the Riemann-Stieltjes
sense. This means that discontinuities of $\boldsymbol{L}(\tau)$
at $\tau_{i}$ represent discrete values of $\frac{\mathrm{d}}{\mathrm{d}t}f_{c}(t-\tau_{i})$.

\section{Contact force calculation}

In what follows, we interpret the meaning of (\ref{eq:redLinForce})
for impact problems. In the simple case when impact occurs with a
stationary stop the contact point must have both a constant position
$y_{1}=\overline{y}_{1}$ and zero velocity $y_{2}=0$. This means
that at the time of initial contact the acceleration of the contact
point becomes infinite as the velocity resets to zero. Despite the
infinite acceleration, the zero mass of the contact point guarantees
a finite contact force. We show that in the transformed model the
zero mass of contact point is equivalent to the condition 
\begin{equation}
\lim_{\tau\to0+}\left[\boldsymbol{L}(\tau)\right]_{2}\neq0.\label{eq:L1cond}
\end{equation}
In other words, if condition (\ref{eq:L1cond}) is satisfied the contact
force is finite. This is a similar, but more general condition that
has been derived by Wang and Kim \cite{WangKim1996}. For multiple
resolved coordinates equation (\ref{eq:L1cond}) must hold for the
coordinates that experience a discontinuity at impact, e.g., the velocities.
To show that our conclusion holds we integrate equation (\ref{eq:redLinForce})
through the initial contact to get 
\[
\boldsymbol{y}(t_{0}^{+})-\boldsymbol{y}(t_{0}^{-})=\boldsymbol{L}^{+}\left(f_{c}(t_{0}^{+})-f_{c}(t_{0}^{-})\right),
\]
where $\boldsymbol{L}^{+}=\lim_{\tau\to0+}\boldsymbol{L}(\tau)$ and
$t_{0}^{-}$ signals a limit from the left and $t_{0}^{+}$ a limit
from the right of the impact. Because the contact force before the
impact $f_{c}(t_{0}^{-})=0$ and the velocity of the contact point
after the impact $y_{2}(t_{0}^{+})=0$, it follows that the initial
contact force is 
\begin{equation}
f_{c}(t_{0}^{+})=-\frac{y_{2}(t_{0}^{-})}{\left[\boldsymbol{L}^{+}\right]_{2}}.\label{eq:firstContactForce}
\end{equation}
During contact $\boldsymbol{y}(t)=\overline{\boldsymbol{y}}=(\overline{y}_{1},0)^{T}$
is constant$ $, which can be substituted into \eqref{eq:redLinForce}
to find out the contact force. Rearranging the resulting equation
yields 
\begin{equation}
\left[\boldsymbol{L}^{+}\right]_{2}\frac{\mathrm{d}}{\mathrm{d}t}f_{c}(t)=-\left[\boldsymbol{L}^{\infty}f_{c}(t)+\boldsymbol{A}\bar{\boldsymbol{y}}+\boldsymbol{g}(t)\right]_{2}-\int_{0+}^{t}\mathrm{d}_{\tau}\left[\boldsymbol{L}(\tau)\right]_{2}\frac{\mathrm{d}}{\mathrm{d}t}\left[f_{c}(t-\tau)\right].\label{eq:intEvolve}
\end{equation}
Equation \eqref{eq:intEvolve} is a delay-differential equation that
contains the history of $f_{c}(t)$, which means that previous impacts
have a great influence on the evolution of the contact force. Clearly,
the contact force is not a continuos function, therefore the derivative
of its history can become infinite. A short calculation shows that
a jump of magnitude $f_{c}^{\mathrm{jump}}$ of $f_{c}$ at $t-\tau^{\star}$
contributes $ $a finite value $\nicefrac{\mathrm{d}}{\mathrm{d}\tau}\boldsymbol{L}(\tau^{\star})f_{c}^{\mathrm{jump}}$
to the integral in (\ref{eq:intEvolve}). 

In conservative systems where shock waves are present, e.g., for an
undamped string, $\boldsymbol{L}(\tau)$ has further isolated discontinuities.
Let $\tau_{d}\neq0$ be the position of such a discontinuity of $\boldsymbol{L}(\tau)$.
If at time $t_{1}=t_{0}+\tau_{d}$ the two bodies are in contact the
contact force develops a further jump. This can be seen by integrating
equation \eqref{eq:intEvolve} for the infinitesimal time interval
$[t_{1}^{-},t_{1}^{+}]$. The integral of the integral on the right
side of \eqref{eq:intEvolve} becomes 
\begin{equation}
\int_{t_{1}^{-}}^{t_{1}^{+}}\int_{0+}^{t}\mathrm{d}_{\tau}\left[\boldsymbol{L}(\tau)\right]_{2}\frac{\mathrm{d}}{\mathrm{d}t}\left[f_{c}(t-\tau)\right]\mathrm{d}t=\left[\int_{0+}^{t}\mathrm{d}_{\tau}\left[\boldsymbol{L}(\tau)\right]_{2}f_{c}(t-\tau)\right]_{t_{1}^{-}}^{t_{1}^{+}},\label{eq:IntOfInt}
\end{equation}
while the other terms are continuous and their integral vanishes.
The right side of equation \eqref{eq:IntOfInt} is regular because
all of its terms are finite. The discontinuity of $\boldsymbol{L}(\tau)$
at $\tau_{d}$ contributes $\left(\boldsymbol{L}(\tau_{d}^{+})-\boldsymbol{L}(\tau_{d}^{-})\right)f_{c}(t-\tau_{d})$
to the integral \eqref{eq:IntOfInt}. Note that $t_{1}-\tau_{d}=t_{0}$,
therefore \eqref{eq:IntOfInt} evaluates to $\left(\boldsymbol{L}(\tau_{d}^{+})-\boldsymbol{L}(\tau_{d}^{-})\right)f_{c}(t_{0}^{+})$
(as $f_{c}(t_{0}^{-})=0$), which means that the contact force at
$t_{1}$ also develops a further discontinuity 
\[
\boldsymbol{L}^{+}\left(f_{c}(t_{1}^{+})-f_{c}(t_{1}^{-})\right)=\left(\boldsymbol{L}(\tau_{d}^{+})-\boldsymbol{L}(\tau_{d}^{-})\right)f_{c}(t_{0}^{+}).
\]

\subsection{Numerical solution of the reduced model\label{sub:numMethod}}

The non-smooth delay-differential equations (\ref{eq:redLinForce},\ref{eq:firstContactForce},\ref{eq:intEvolve})
are somewhat unusual and therefore standard numerical techniques are
not directly applicable to them. In what follows, we use a simple
explicit Euler method and the rectangle rule of numerical integration
to approximate the solution of \eqref{eq:redLinForce} for cases when
$\left[\boldsymbol{L}^{+}\right]_{2}\neq0$. We assume that time is
quantized in $\varepsilon$ chunks, so that $\boldsymbol{y}_{q}=\boldsymbol{y}(q\varepsilon)$,
$f_{c,q}=f_{c}(q\varepsilon)$, where $q=0,1,2,\ldots$. If there
is no contact and consequently no contact force ($f_{c,q}=0$), the
only unknown is the state variable $\boldsymbol{y}_{q}$. Therefore
the evolution of the resolved variables is given by
\begin{equation}
\boldsymbol{y}_{q+1}=\boldsymbol{y}_{q}+\varepsilon\left(\boldsymbol{A}\boldsymbol{y}_{q}+\boldsymbol{g}(q\varepsilon)\right)+\sum_{j=0}^{q-1}\left(\boldsymbol{L}_{j+1}-\boldsymbol{L}_{j}\right)\left(f_{c,q-j}-f_{c,q-j-1}\right),\label{eq:NumStateUpdate}
\end{equation}
where $\boldsymbol{L}_{j}=\boldsymbol{L}(j\varepsilon)$ and $\boldsymbol{L}_{0}=\boldsymbol{L}^{+}$.
Contact of the impacting bodies is detected, when $\left[\boldsymbol{y}_{q+1}\right]_{1}\le\overline{y}_{1}$.
In this case the resolved variables are kept constant with $\boldsymbol{y}_{q+1}=\overline{\boldsymbol{y}}$.
Also, equation \eqref{eq:firstContactForce} is applied at the onset
of contact, so that the initial contact force becomes 
\begin{equation}
f_{c,q+1}=-\frac{\left[\boldsymbol{y}_{q}\right]_{2}}{\left[\boldsymbol{L}^{+}\right]_{2}}.\label{eq:NumFirstContact}
\end{equation}
The subsequent values of the contact force are calculated by 
\begin{equation}
f_{c,q+1}=f_{c,q}-\frac{\varepsilon}{\left[\boldsymbol{L}^{+}\right]_{2}}\left[\boldsymbol{L}^{\infty}f_{c,q}+\boldsymbol{A}\bar{\boldsymbol{y}}+\boldsymbol{g}(q\varepsilon)\right]_{2}-\frac{1}{\left[\boldsymbol{L}^{+}\right]_{2}}\sum_{j=0}^{q-1}\left[\boldsymbol{L}_{j+1}-\boldsymbol{L}_{j}\right]_{2}\left(f_{c,q-j}-f_{c,q-j-1}\right).\label{eq:FcUpdate}
\end{equation}
which is the discretized counterpart of \eqref{eq:intEvolve}. If
$f_{c,q+1}$ as predicted by equation \eqref{eq:FcUpdate} becomes
negative, we set $f_{c,q+1}=0$ and continue the calculation with
\eqref{eq:NumStateUpdate}.

\section{Impacting cantilever beam models}

\begin{figure}[t]
\begin{centering}
\includegraphics[scale=1.2]{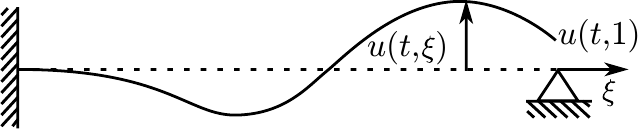} 
\par\end{centering}

\caption{\label{fig:beam}Impacting cantilever beam.}
\end{figure}

Our theory sets a criterion in the form of equation (\ref{eq:L1cond})
for the regularity of the mechanical model, which can be used to test
different models of elastic structures. We consider the example of
a cantilever beam in Fig.~\ref{fig:beam} described by two different
models. Through our calculation it becomes clear why the Euler-Bernoulli
model often used in impact models \cite{Melcher,Fathi1994,YinQinZou2007}
exhibits signs of singularity.

\begin{figure}
\begin{centering}
\includegraphics[width=0.6\linewidth]{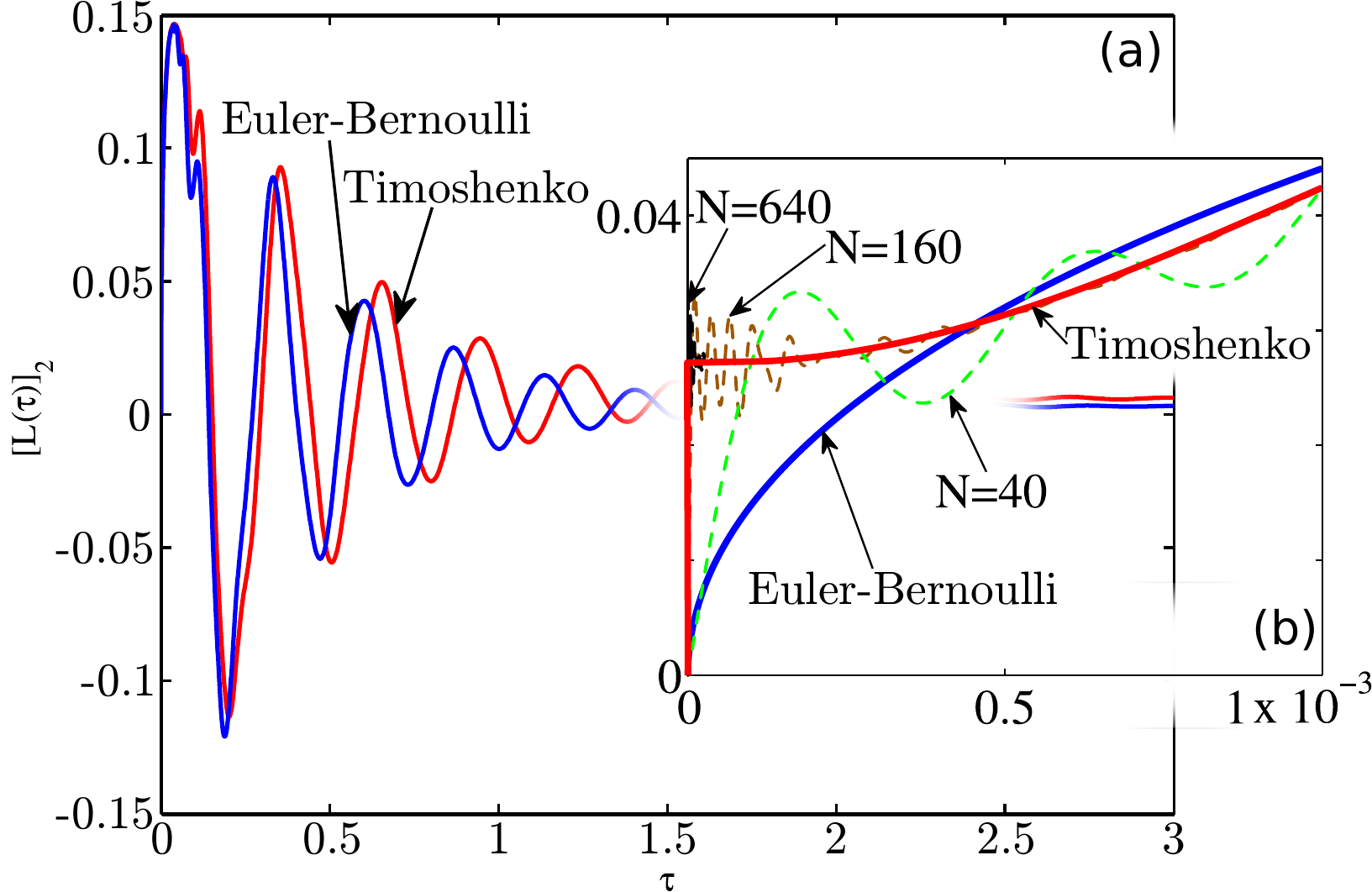} 
\par\end{centering}

\caption{(color online) Graph of $\left[\boldsymbol{L}(\tau)\right]_{2}$ for
the Euler-Bernoulli and the Timoshenko beam models. The inset shows
that the Timoshenko beam model converges to a function where $\left[\boldsymbol{L}^{+}\right]_{2}\neq0$,
while the Euler-Bernoulli model is singular since its $\boldsymbol{L}(\tau)$
is continuos. The parameters are $\beta=4800$, $\gamma=\nicefrac{1}{4}$.\label{fig:L1fun}}
\end{figure}

As the underlying elastic structure, consider the Euler-Bernoulli
cantilever beam model 
\begin{equation}
\frac{\partial^{2}u}{\partial t^{2}}=-\frac{\mathrm{\partial}^{4}u}{\partial\xi^{4}}+\psi_{2}f_{e}(t),\; u(t,0)=\left.\frac{\mathrm{\partial}u(t,\xi)}{\partial\xi}\right|_{\xi=0}=\left.\frac{\mathrm{\partial}^{2}u(t,\xi)}{\partial\xi^{2}}\right|_{\xi=1}=0,\label{eq:EBbeam}
\end{equation}
with $\left.\nicefrac{\mathrm{\partial}^{3}u(t,\xi)}{\partial\xi^{3}}\right|_{\xi=1}=f_{c}(t)$,
where $u(t,\xi)$ represents the deflection of the beam. The natural
frequencies of (\ref{eq:EBbeam}) are determined by the equation $1+\cos\sqrt{\omega_{k}}\cosh\sqrt{\omega_{k}}=0$,
while the mode shape values at the end of the beam are described by
\begin{equation}
\boldsymbol{n}=(2,-2,2,-2,\ldots)^{T}.\label{eq:nVecEB}
\end{equation}
 On the other hand the Timoshenko beam model is represented by 
\begin{align}
\frac{\partial^{2}u}{\partial t^{2}} & =\beta\gamma\left(\frac{\partial^{2}u}{\partial\xi^{2}}-\frac{\partial\phi}{\partial\xi}\right)+\psi_{2}f_{e}(t),\nonumber \\
\frac{\partial^{2}\phi}{\partial t^{2}} & =\beta\frac{\partial^{2}\phi}{\partial\xi^{2}}+\beta^{2}\gamma\left(\frac{\partial u}{\partial\xi}-\phi\right),\nonumber \\
u(t,0) & =\phi(t,0)=\phi(t,1)=0,\label{eq:TMbeam}
\end{align}
and $\left.\nicefrac{\mathrm{\partial}}{\partial\xi}u(t,1)\right|_{\xi=1}-\phi(t,1)=f_{c}(t)$,
where $f_{e}(t)$ is an external forcing through the second mode shape
$\psi_{2}$, and $\phi$ is the rotation angle of the cross-section
of the beam. The resolved coordinates in both cases are $y_{1}=u(t,1)$
and $y_{2}=\nicefrac{\mathrm{\partial}}{\partial t}u(t,1)$. In case
of the Timoshenko beam, modal decomposition in the form of (\ref{eq:modeDecomp})
is not practical, instead we use Chebyshev collocation \cite{trefethen}
to discretize the system with $N$ number of collocation points. This
is possible, since our formulation is not restricted to modal decomposition,
the matrix $\boldsymbol{R}$ can represent any form of discretization.
The governing equations (\ref{eq:EBbeam},\ref{eq:TMbeam}) are conservative,
thus to obtain a decaying solution we add modal damping ratios $D_{k}=\nicefrac{1}{10}$
to both systems (\ref{eq:EBbeam},\ref{eq:TMbeam}) after being discretized.
The result of our calculation is shown in Fig.~\ref{fig:L1fun}.
It can be seen that for the Euler-Bernoulli model $\left[\boldsymbol{L}(\tau)\right]_{2}$
is continuos while for the Timoshenko model $\left[\boldsymbol{L}(\tau)\right]_{2}$
is discontinuous. This result is in accordance with the fact that
the scaling exponent of the highest natural frequencies for the Euler-Bernoulli
model is $\alpha=2$ and for the Timoshenko model $\alpha=1$ \cite{Rensburg2006}.
This means that the Euler-Bernoulli model is not suitable for impact
calculations.

\begin{figure}
\begin{centering}
\includegraphics[width=0.99\linewidth]{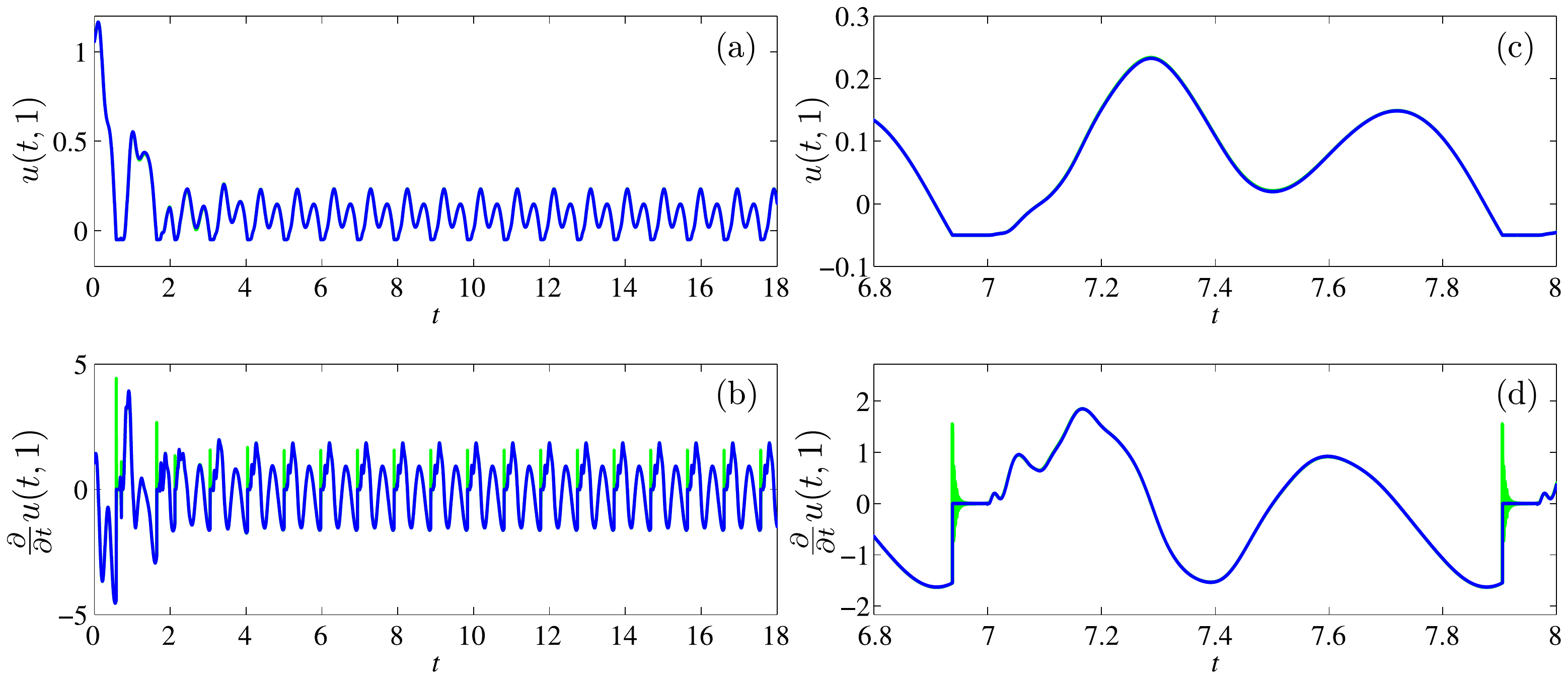} 
\par\end{centering}

\caption{(color online) Vibrations of the impacting Timoshenko beam using two
impact models. The rigid stop as illustrated in Fig.~\ref{fig:beam}
is placed at $\overline{y}_{1}=-0.05$, and the forcing is $f_{e}(t)=30\cos\left(13t\right)$
through the second mode. Trajectories of the reduced model (\ref{eq:redLinForce},\ref{eq:firstContactForce},\ref{eq:intEvolve})
are shown in dark (blue) and the solution of the CoR model (\ref{eq:modeDecomp},\ref{eq:CoRvelocity})
is represented by light (green) lines for comparison. The time step
used to solve the reduced model is $\varepsilon=3.5\times10^{-5}$
and the number of collocation points used to solve the CoR model is
$N=20$. Panels (a,b) show that the solution converges to a periodic
orbit. A single period of the solution is illustrated in panels (c,d).\label{fig:impsol}}
\end{figure}

\begin{figure}
\begin{centering}
\includegraphics[width=0.6\linewidth]{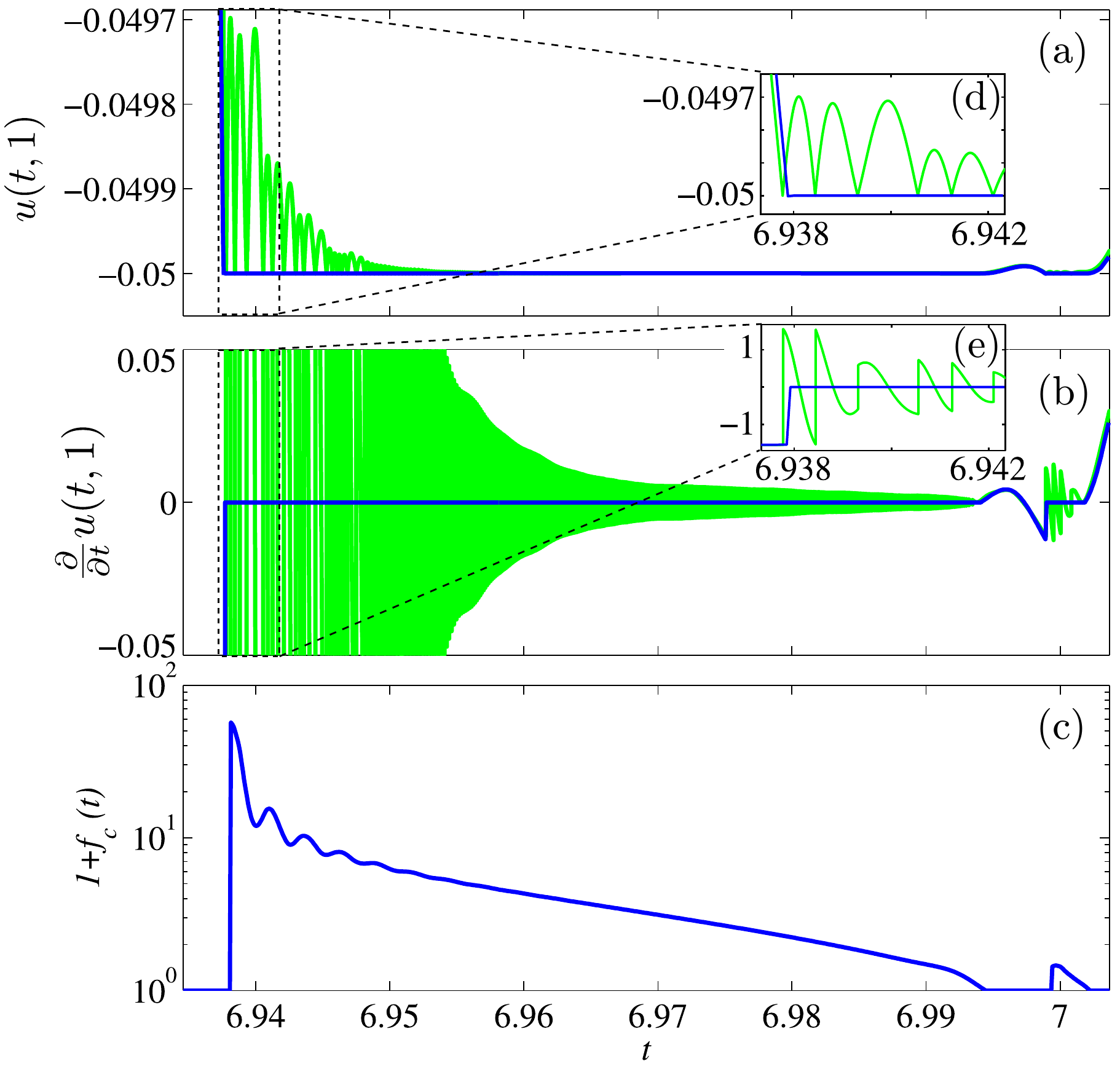}
\par\end{centering}

\caption{(color online) A sequence of two impacts within a period of the periodic
solution in Fig. \ref{fig:impsol}. The solution of the CoR model
(\ref{eq:modeDecomp},\ref{eq:CoRvelocity}) as shown by the light
(green) lines is highly oscillatory. The dark (blue) lines show that
the reduced model (\ref{eq:redLinForce},\ref{eq:firstContactForce},\ref{eq:intEvolve})
is similarly accurate and eliminates the high-frequency chatter of
the CoR model. (a,b). Insets (d,e) show the small-scale dynamics of
the CoR model. The contact force is finite and continuos after the
initial contact (c).\label{fig:whileImpact}}

\end{figure}

\subsection{The Timoshenko model}

In the previous section we have shown that the Timoshenko beam model
\eqref{eq:FcUpdate} satisfies our criterion \eqref{eq:L1cond}, which
guarantees that equations (\ref{eq:redLinForce},\ref{eq:firstContactForce},\ref{eq:intEvolve})
are well defined and that the contact force is finite during impact.
To confirm that this is indeed the case we simulate the impacting
cantilever beam with the Timoshenko beam model \eqref{eq:TMbeam}
as illustrated in figure \ref{fig:beam}. The rigid stop is placed
at $\overline{y}_{1}=-0.05$, so that the beam contacts the stop when
the position of its tip reaches $u(t,1)=-0.05$. In all our simulations
we use the initial conditions $u(0,\xi)=\nicefrac{\mathrm{\partial}}{\partial t}\left|u(t,\xi)\right|_{t=0}=1.056\psi_{1}(\xi)$,
where $\psi_{1}$ is the first mode shape of the structure and it
is normalized by $\psi_{1}(1)=1$. We also force the beam through
its second mode with $f_{e}(t)=30\cos\left(13t\right)$. The numerical
solution is obtained using equations (\ref{eq:NumStateUpdate},\ref{eq:NumFirstContact},\ref{eq:FcUpdate}).
In addition to our method we also simulate the dynamics using the
CoR model described in \cite{Vyasarayani2010}. This comparison highlights
that high-frequency chatter is eliminated when our method is used.

The CoR model that we are using for comparison describes the impact
as an infinitesimally short process. During the impact an impulse
is applied at the contact point that alters the velocity state of
the body. The magnitude of this impulse is determined by the desired
rebound velocity of the contact point which is $-C_{R}$ times the
incident velocity. Assuming that the equation of motion is in the
form of \eqref{eq:modeDecomp} and impact occurs at $t_{0}$, the
after-impact velocity of the structure is
\begin{equation}
\dot{\boldsymbol{x}}(t_{0}^{+})=\left(\boldsymbol{I}-(1+C_{R})\frac{\boldsymbol{n}\otimes\boldsymbol{n}}{\boldsymbol{n}^{2}}\right)\dot{\boldsymbol{x}}(t_{0}^{-}),\label{eq:CoRvelocity}
\end{equation}
where $\otimes$ means the outer product between vectors. The position
of the structure remains the same throughout the impact: $\boldsymbol{x}(t_{0}^{+})=\boldsymbol{x}(t_{0}^{-})$.
The CoR model is physically questionable because it treats the impact
as an infinitesimally short event. It can however, reproduce most
experimental observations \cite{Popp}. The source of high-frequency
chatter can be explained by equation \eqref{eq:CoRvelocity}, which
stipulates that the change in modal velocities is proportional to
a constant times vector $\boldsymbol{n}$. Elements of $\boldsymbol{n}$
corresponding to high frequency modes are of the same magnitude as
for low frequency modes (e.g., see equation \eqref{eq:nVecEB}), which
means that after an impact the tip of the beam acquires a high frequency
vibration. Due to this vibration another impact is likely to follow
shortly and repeatedly, which results in chatter that has roughly
the same frequency as the highest vibration mode of the structure.
This is illustrated by the light (green) lines in figure \ref{fig:whileImpact}.

In contrast to the CoR method (\ref{eq:modeDecomp},\ref{eq:CoRvelocity})
our method as solved by equations (\ref{eq:NumStateUpdate},\ref{eq:NumFirstContact},\ref{eq:FcUpdate})
eliminates chatter and produces a much smoother result, which are
shown by dark (blue) lines Figs.~\ref{fig:impsol} and \ref{fig:whileImpact}.
On the larger scale in figure \ref{fig:impsol} the two solutions
roughly coincide, while on the smaller scale in figure \ref{fig:whileImpact}
the high frequency chatter is apparent and would increase in frequency
if more vibration modes or collocation points on the beam were used.
This high frequency chatter can stall numerical simulations, while
the time-step in our method is not affected by the inclusion of higher
natural frequencies.

The contact force in the CoR model is infinite at times of contact
and hence it cannot be calculated. On the other hand our method allows
the calculation of the contact force, which is finite as shown in
figure \ref{fig:whileImpact}(c). The contact force even becomes a
smooth function of time after the onset of contact.

%The contact force as a function of the number of modes  can be seen in figure . As  increases, the initial jump of the contact force decreases as  develops its discontinuity at  Also notice that the time variation of the contact force gets smoother.

\section{Conclusions}

In this paper we introduced a new way of modeling the impact mechanics
of elastic structures. With our method regularity of the model can
be predicted and a finite and piecewise continuos contact force can
be calculated. The key to this result is that the delay equation description
preserves the infinite dimensional nature of the mechanics and the
zero mass of the contact point. The results presented in this paper
open a significant number of new avenues of research. Models that
show non-deterministic behavior such as the Painleve paradox \cite{PainleveAlan}
might be regularized through our method. The strong dependence of
dynamical phenomena on the number of underlying dimensions \cite{diBernardoBook}
could also be eliminated, since our framework considers all the infinite
dimensions. Further, the bifurcation theory of non-smooth delay equations
requires attention in order to understand the implications of our
regularized impact mechanics, especially in how far it is an improvement
over finite dimensional models.
\begin{acknowledgments}
The author thanks Gábor Stépán, who brought his attention to the work
of Chorin et al. \cite{ChorinPNAS}. He also thanks Alan R. Champneys,
John Hogan and Thibaut Putelat for useful discussion and comments
on the manuscript. Corrections to the text by Galit Szalai are greatly
appreciated. 
\end{acknowledgments}
\appendix
%dummy comment inserted by tex2lyx to ensure that this paragraph is not empty

\section{Contact force asymptotics}

In this appendix we approximate the contact force at the onset of
an impact. We assume that an impact takes place at $t=t_{0}$. To
help the notation we define $\square^{-}=\square(t_{0})$ and $\square^{+}=\square(t_{0}+\delta t)$,
where $\square$ stands for any dependent variable. We aim to calculate
a constant contact force $f_{c}$ that allows the elastic body to
overlap with the rigid stop for an exactly $\delta t$ long time interval.
As $\delta t$ tends to zero the overlap is removed, hence the calculated
$f_{c}$ force tends to the actual contact force.

\begin{figure}[t]
\begin{centering}
\includegraphics[width=0.49\linewidth]{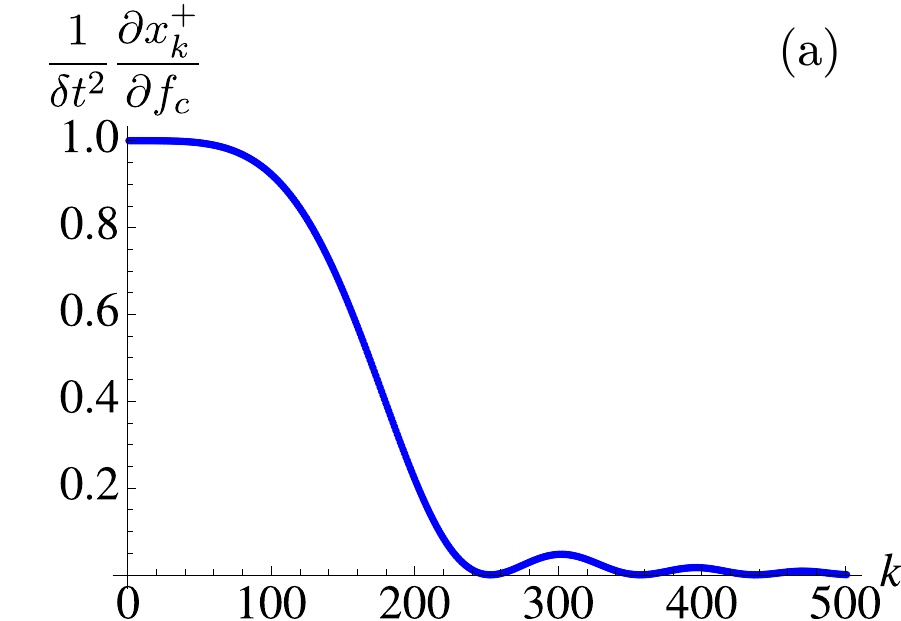}\includegraphics[width=0.49\linewidth]{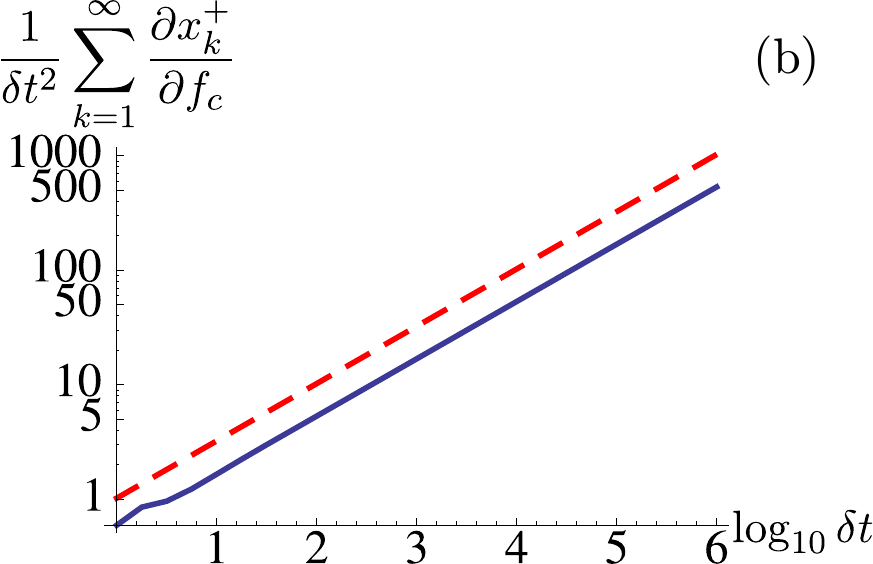} 
\par\end{centering}

\caption{\label{fig:coeffFig}(a) Coefficients (\ref{eq:DfXoft}) for $\delta t=10^{-5}$,
$D_{k}=0$ and $\omega_{k}=\left(k\pi-\frac{\pi}{2}\right)^{2}$.
(b) The sum of the infinite series of coefficients (\ref{eq:DfXoft})
is shown by the continuos (blue) line. The dashed (red) line is the
estimate of $N$, in Eqn. (\ref{eq:Napprox}) that overestimates the
sum.}
\end{figure}

To calculate this constant $f_{c}$ one needs to solve 
\[
0=\boldsymbol{n}\cdot\boldsymbol{x}^{+}
\]
for $f_{c}$. Due to the linearity of equation (\ref{eq:modeDecomp}),
the motion depends linearly on the contact force. Therefore expanding
the constraint $0=\boldsymbol{n}\cdot\boldsymbol{x}^{+}$ at $f_{c}=0$
we get the exact equation 
\begin{equation}
0=\boldsymbol{n}\cdot\boldsymbol{x}_{f_{c}=0}^{+}+f_{c}\sum_{k=1}^{N}\psi_{k}(\chi^{\star})\frac{\partial x_{k}^{+}}{\partial f_{c}}.\label{eq:fcCondition}
\end{equation}
The derivatives in (\ref{eq:fcCondition}) are calculated in closed
form as 
\begin{equation}
\frac{\partial x_{k}^{+}}{\partial f_{c}}=-\frac{\psi_{k}(\chi^{\star})}{\omega_{k}^{2}}\Biggl(\frac{\mathrm{e}^{-D_{k}\omega_{k}\delta t}}{\sqrt{1-D_{k}^{2}}}\sin\left(\omega_{k}\sqrt{1-D_{k}^{2}}\delta t\right)+\mathrm{e}^{-D_{k}\omega_{k}\delta t}\cos\left(\omega_{k}\sqrt{1-D_{k}^{2}}\delta t\right)-1\Biggr).\label{eq:DfXoft}
\end{equation}
When calculating the derivatives (\ref{eq:DfXoft}) for $\delta t=10^{-5}$,
we get a vanishing sequence as is illustrated in Fig. \ref{fig:coeffFig}(a).

We assume that before the impact, the structure has a smooth motion
so that the displacement without the contact force can be approximated
as 
\[
\boldsymbol{n}\cdot\boldsymbol{x}_{f_{c}=0}^{+}\approx\delta t\,\boldsymbol{n}\cdot\dot{\boldsymbol{x}}^{-}.
\]
To quantify how the derivatives in Eqn. (\ref{eq:DfXoft}) vanish
we asymptotically expanded them as 
\begin{align*}
\frac{\partial x_{k}^{+}}{\partial f_{c}} & =\frac{1}{2}\psi_{k}(\chi^{\star})\delta t^{2}\left(1-\frac{2}{3}D_{k}\omega_{k}\delta t-\frac{1-4D_{k}^{2}}{12}\omega_{k}^{2}\delta t^{2}+\cdots\right),
\end{align*}
Substituting these two estimates into (\ref{eq:fcCondition}) we get
the contact force 
\begin{equation}
f_{c}\approx\frac{-2\boldsymbol{n}\cdot\dot{\boldsymbol{x}}^{-}}{\delta t\sum_{k=1}^{N}\psi_{k}^{2}(\chi^{\star})}=\frac{-2\boldsymbol{n}\cdot\dot{\boldsymbol{x}}^{-}}{N\delta t\overline{\psi_{k}^{2}(\chi^{\star})}},\label{eq:fcform}
\end{equation}
up to the leading order, where the limit $N$ of the summation equals
the smallest $k$ for which (\ref{eq:DfXoft}) vanishes, and $\overline{\psi_{k}^{2}(\chi^{\star})}$
is the average value of $\psi_{k}^{2}(\chi^{\star})$, for $k\le N$.
$N$ can be calculated as the smallest $k$ such that the expanded
coefficients become small 
\[
\frac{1}{\delta t^{2}}\frac{\partial x_{k}^{+}}{\partial f_{c}}\approx1-\frac{2}{3}D_{k}\omega_{k}\delta t-\frac{1-4D_{k}^{2}}{12}\omega_{k}^{2}\delta t^{2}<\eta\ll1.
\]
Using the asymptotic scaling of the natural frequencies $\omega_{k}=\omega_{0}k^{\alpha}$
and assuming zero damping ($D_{k}=0$) we get 
\begin{equation}
N>\left(\frac{2\sqrt{3-3\eta}}{\omega_{0}}\right)^{\frac{1}{\alpha}}\delta t^{-\frac{1}{\alpha}},\label{eq:Napprox}
\end{equation}
therefore the restoring force is 
\[
f_{c}=-C\delta t^{\frac{1}{\alpha}-1}\boldsymbol{n}\cdot\dot{\boldsymbol{x}}^{-},
\]
where $C\approx2\left(\frac{2\sqrt{3-3\eta}}{\omega_{0}}\right)^{-\frac{1}{\alpha}}\overline{\psi_{k}^{2}(\chi^{\star})}^{-1}.$
In order to check validity of our estimate of $N$ we plotted the
sum of derivatives in (\ref{eq:DfXoft}) and compared to the estimate
(\ref{eq:Napprox}) in Fig. \ref{fig:coeffFig}(b).

As the last step we calculate the change in modal velocities 
\begin{align*}
\dot{x}_{k}^{+} & =\dot{x}_{k}^{-}-C\delta t^{\frac{1}{\alpha}-1}\boldsymbol{n}\cdot\dot{\boldsymbol{x}}^{-}\frac{\partial\dot{x}_{k}^{+}}{\partial f_{c}}\\
 & \approx\dot{x}_{k}^{-}-\psi_{k}(\chi^{\star})C\delta t^{\frac{1}{\alpha}}\boldsymbol{n}\cdot\dot{\boldsymbol{x}}^{-},
\end{align*}
which means that if $\delta t\to0$, there is no change in individual
modal velocities, except when $\alpha\to\infty$, that is the case
of a rigid body. However, when calculating the velocity of the impacting
point after the impact we have 
\[
\boldsymbol{n}\cdot\dot{\boldsymbol{x}}^{+}=-\boldsymbol{n}\cdot\dot{\boldsymbol{x}}^{-}.
\]

\bibliographystyle{plain}
\bibliography{AllRef}

\end{document}